\documentclass[a4paper,12pt]{article}
\usepackage[english]{babel}
\usepackage{amsmath,amsfonts,amssymb,amscd,color,epsfig,amsthm,enumerate,calc}

\newcommand{\diam}{\ensuremath{\operatorname{diam}}}
\newcommand{\D}{\ensuremath{{\mathbb   D}}}

\newcommand{\eps}{{\epsilon}}
\newcommand{\epsk}{{\epsilon^k}}
\newcommand{\epsn}{{\epsilon^n}}
\newcommand{\ga}{\ensuremath{{\gamma}}}

\newcommand{\itemref}[1]{\ref{#1}.}
\newcommand{\id}{\ensuremath{{\operatorname{id}}}}
\newcommand{\mapfromto}[3]{\hbox{\ensuremath{#1 : #2 \longrightarrow #3}}}

\newcommand{\N}{\ensuremath{\mathbb  N}}

\newcommand{\Sm}{\ensuremath{{\setminus}}}

\newcommand{\Sig}{\ensuremath{{\Sigma}}}
\newcommand{\Sigd}{\ensuremath{{\Sigma_d}}}
\newcommand{\Sigdk}{\ensuremath{{\Sigma_d^k}}}

\newcommand{\veps}{\ensuremath{{\underline{\epsilon}}}}

\newcommand{\wheps}{{\widehat{\epsilon}}}
\newcommand{\whepsk}{{\widehat{\epsilon}^k}}
\newcommand{\whga}{\ensuremath{\widehat\ga}}
\newcommand{\whveps}{{\ensuremath{\widehat{\underline{\epsilon}}}}}
\newcommand{\whV}{\ensuremath{\widehat V}}
\def\ga{\gamma}

\def\eps{\epsilon}

\newtheorem{theorem}{Theorem}[section]
\newtheorem{lemma}[theorem]{Lemma}
\newtheorem{proposition}[theorem]{Proposition}
\newtheorem{corollary}[theorem]{Corollary}

\newtheorem{definition}[theorem]{Definition}

\newcommand{\THM}{\begin{theorem}}
\newcommand{\ENUM}{\begin{enumerate}}
\newcommand{\ENDENUM}{\end{enumerate}}
\newcommand{\REFTHM}[1] { \begin{theorem}\label{#1} }
\newcommand{\RREFTHM}[2] { \begin{theorem}[#1]\label{#2} }
\newcommand{\ENDTHM}{\end{theorem}}
\newcommand{\REFPROP}[1]{\begin{proposition}\label{#1} }
\newcommand{\RREFPROP}[2]{\begin{proposition}[#1]\label{#2} }
\newcommand{\PROP}{\begin{proposition}}
\newcommand{\ENDPROP}{\end{proposition} }
\newcommand{\REFDEF}[1]{\begin{definition}\label{#1} }
\newcommand{\DEF}{\begin{definition}}
\newcommand{\ENDDEF}{\end{definition} }
\newcommand{\REFLEM}[1]{\begin{lemma}\label{#1} }
\newcommand{\RREFLEM}[2]{\begin{lemma}[#1]\label{#2} }
\newcommand{\LEM}{\begin{lemma}}
\newcommand{\ENDLEM}{\end{lemma} }
\newcommand{\REFCOR}[1]{\begin{corollary}\label{#1} }
\newcommand{\COR}{\begin{corollary}}
\newcommand{\ENDCOR}{\end{corollary} }

\newcommand{\thmref}[1]{Theorem~\ref{#1}}
\newcommand{\propref}[1]{Proposition~\ref{#1}}

\newcommand{\PROOF}{\begin{proof}}
\newcommand{\ENDPROOF}{\end{proof}}

\begin{document}
\title{Polynomial-like semi-conjugates\\
of the shift map.
\footnote{This work was supported by grant 272-07-0321 from
the Research Council for Nature and Universe and by the Marie Curie RTN CODY Project No 35651.}}
\author{Carsten Lunde Petersen}

\maketitle
\begin{abstract}
In this paper I prove that for a polynomial of degree $d$ with a Cantor Julia set $J$, the
Julia set can be understood as the simplest possible quotiont of the one sided shift
space $\Sigma_d$ with dynamics given by the shift. 
Here simplest possible means that, 
the projection {\mapfromto \pi {\Sigma_d} J} is as injective as possible.
\end{abstract}

\section{Introduction}
Denote by $\Sigd = {\{0, \ldots , d-1\}}^{\N}$ the set of one-sided infinite sequences 
of symbols the $0, \ldots , d-1$ equipped with the natural product topology. 
And denote by {\mapfromto \sigma {\Sigma_d} {\Sigma_d}} the
shift map:
$$
\sigma({(\eps_i)}_i) = {(\eps_{i+1})}_i = (\eps_2,\eps_3, \ldots ).
$$

Douady and Hubbard introduced in \cite{DouadyHubbard} the notion of polynomial-like maps. 
Here we shall use a slightly generalized version of such maps (see also \cite{LyubichVolberg}):

Let {\mapfromto f {U'} U} be a proper holomorphic map where $U\simeq \D$, 
$U'\subset\subset U$, $U'= U_1\cup U_1\cup\ldots\cup U_N$, $U_i\simeq \D$ for 
each $i$ and $U_i\cap U_j = \emptyset$ for $i\not=j$.

The filled-in Julia set $K_f$ for $f$ is the set of points:
$$
K_f = \{z\in U'| f^n(z)\in U', \forall\; n\in\N\},
$$
and the Julia set is its topological boundary $J_f=\partial K_f$.

For such a map the degree $d$ is the sum of the degrees $d_i$ of the restrictions
{\mapfromto {f_|} {U_i} U}. By the Riemann-Hurwitz formula $f$ has counting multiplicity $d_i'=d_i-1$
critical points in $U_i$. In particular if $f$ does not have any critical point in some
$U_i$, then $f$ has a globally defined inverse branch {\mapfromto {f_i=f_|^{-1}} U {U_i}}. 
In particular if $f$ has no critical points at all then $d=N$ and $f$ has $d$ distinct
globally defined inverse branches. 
In this case it follows that $K_f=J_f$ is a Cantor set and an
elementary proof going back to Fatou shows that in the later case there is a
homeomorphism {\mapfromto \pi {\Sigma_d} {J_f}} such that $\pi\circ\sigma=f\circ\pi$.

If no critical point of $f$ is periodic then the function {\mapfromto \chi {K_f} \N} given
by the maximal local degree of iterates of $f$ near $z$:
$$
\chi(z) = \sup_{n\in\N}\deg(f^n,z)
$$
is bounded by the product of the local degrees of $f$ at its critical points and satisfies 
$$
\chi(z) = \deg(f,z)\cdot\chi(f(z)).
$$
Since $f$ has $d'=d-N = \sum d_i'$ critical points the function $\chi$ is bounded by $2^{d'}$.
Note that any periodic critical point is surrounded by an open attracted basin,
and thus belongs to the interior of $K_f$.

The main theorem of this paper is:

\REFTHM{maintheorem} 
Let {\mapfromto f {U'} U} be a degree $d>1$ generalized polynomial-like map in the sense above. 
If $K_f=J_f$ is a Cantor set containing all critical points of $f$. 
Then there is a semi-conjugacy 
{\mapfromto \pi {\Sigma_d} {J_f}}, $\pi\circ\sigma = f\circ \pi$ such that:
$$
\forall z\in J_f : \#\pi^{-1}(z) = \chi(z).
$$
\ENDTHM 
\noindent\textbf{Remark 1)} Branner and Hubbard proved in \cite{BrannerHubbard} 
that there are many cubic polynomials with a generalized polynomial-like restriction as 
above satisfying the hypothesis and thus the conclusion of the above Theorem. 
Moreover recently this Branner-Hubbard Theorem has been extend to all 
degrees and all orders of critical points. See e.g. \cite{TanYin}, \cite{QuiYin} and \cite{KozlowskiStrien}.\\
\noindent\textbf{Remark 2)} 
The main theorem is related to the structure of the complement of the cubic connectedness locus 
through ther paper \cite{BlanchardDevaneyKeen} of Devaney, Goldberg and Keen.\\
\noindent\textbf{Remark 3)} 
The injectivity statement of the main Theorem is best possible, 
since if for some point $z$ :$\#\pi^{-1}(P(z)) = l$ 
and if the local degree of $f$ at $z$ is $m\geq 1$. Then $\#\pi^{-1}(z) = ml$.\\
\noindent\textbf{Remark 4)} 
$N\geq 2$ since if not $J_f$ would be connected and not a Cantor set.\\
\noindent\textbf{Remark 5)} 
The hypothesis that $J_f$ is Cantor set is equivalent to asking that the diameters of the
connected components of $f^{-n}(U)$ converge to zero as $n$ tends to infinity.

Towards a proof of the main theorem we introduce some notation.\\
\noindent\emph{We shall in the following tacitly assume the hypothesies of the Theorem, i.e. 
{\mapfromto f {U'} U} is a generalized polynomial-like map for which $K_f=J_f$ is a 
Cantor set containing all critical points of $f$}.\\
\noindent 
Note that taking a restriction with $U$ slightly smaller if necessary we can assume the boundaries of all
disks $U$ and $U_i$ are smooth and disjoint. Let $w\in U\Sm U'$ be arbitrary and let
$w_0, \ldots ,w_{d-1}$ denote the $d$ distinct preimages of $w$, and let $i=i(j)$ denote
the function given by $w_j\in U_{i(j)}$. Renumbering if necessary we can assume that $i$
is weakly increasing, i.e. we have filled-in from below.

Let {\mapfromto \phi \D {U\Sm J_f}} be a universal covering with $\phi(0) = w$.

\REFPROP{liftedinverses} There exist $d$ (univalent) lifts {\mapfromto {g_i} \D \D}, 
$i = 0, \ldots , d-1$ of $\phi$ to $f\circ\phi$, i.e.~$f\circ\phi\circ g_i = \phi$ with
$\phi\circ g_i(0) = w_i$. 
These satisfy $\phi\circ g_j(z)\not=\phi\circ g_{j'}(z)$ for $j\not= j' \mod d$, 
i.e.~for any $z\in\D$ the points $\phi\circ g_j(z)$ are the $d$ distinct preimages of $\phi(z)$ under $f$.
In particular 
$$
f^{-1}(\phi(\D)) = f^{-1}(U\Sm J_f) = U'\Sm J_f = \bigcup_{i=0}^{d-1} (\phi\circ g_i)(\D).
$$
\ENDPROP
\noindent Remark that the $g_i$ are by no means unique.

\begin{proof}{}
For $0\leq j <d$ and $i=i(j)$ let $V_i$ be a connected component of $\phi^{-1}(U_i\Sm J_f)$. 
Then $V_i$ is simply connected, because $U_i\Sm J_f$ is a retract of $U\Sm J_f$. 
Hence the restriction {\mapfromto \phi {V_i} {U_i\Sm J_f}} is a universal covering map. 
Since the restriction {\mapfromto f {U_i\Sm J_f} {U\Sm J_f}} 
has no critical points it is also a covering and thus each $f\circ\phi_{|V_i}$ is a universal covering. 
Let $x_j\in V_i$ be any point with $\phi(x_j) = w_j$. 
Then there is a unique lift {\mapfromto {g_j} \D {V_j}} of the universal covering $\phi$ 
to the (universal) covering $f\circ\phi_{|V_j}$ mapping $0$ to $x_j$. 
Being lifts of $f\circ\phi$ to $\phi$, any two of the $g_j$ either egree everywhere or nowhere. 
They are chosen to disagree at $0$.
\ENDPROOF
Note that changing the choice of some $z_j$ to some other preimage $z_j'$ of $w_j$ amounts to post composing 
$g_j$ with the decktransformation  for $\phi$, which maps $z_j$ to $z_j'$. 
We shall think and speak of the maps $g_j$ as lifts of $f^{-1}$ though technically 
they are self-maps of a different space.

For $k\geq 1$ let $\Sigma_d^k$ denote the set of $k$-blocks 
$\eps^k = (\eps_1, \ldots, \eps_k)$ in the alphabet $\{0, \ldots, d-1\}$. Every such $\eps^k$ defines 
a ``cylinder'' clopen set 
$$
\{{(\tau_j)}_j| \tau_j=\eps_j,  j\leq k\}.
$$
The map $\sigma$ thus has a natural extension as a map from $\Sigma_d^k$ to $\Sigma_d^{k-1}$. 
Also for $n<k$ there is a natural projection from $\Sigdk$ to $\Sig_d^n$: $\epsk\to\eps^n$, 
which simply forgets the last $k-n$ entries.

The obvious idea for proving \thmref{maintheorem} 
would now be to iterate the $d$ branches {\mapfromto {g_j} \D \D} of the inverse of $f$,  
project back to $U$ and obtain sets for defining a semiconjugacy.
More precisely for $\veps = (\eps_1, \eps_2, \ldots, \eps_k, \ldots) \in \Sigma_d$ and $k\in\N$ define
$$
g_{\epsk} = g_{\eps_1}\circ \ldots \circ g_{\eps_k},
$$
and
$$
V_{\epsk} = g_{\epsk}(\D).
$$
Then for each $\veps\in\Sigd$ the set 
$\cap_{k\geq 1} \phi(V_{\epsk})$ 
is a connected subset of the Cantor set $J_f$ and thus a singleton $\{z_\veps\}$. 
Define {\mapfromto \Psi {\Sigma_d} {J_f}} by $\Psi(\veps) = z_\veps$. 
Then $\Psi$ is indeed a semi-conjugacy of {\mapfromto \sigma {\Sigd}{\Sigd}} to {\mapfromto P {J_f} {J_f}}. 
However in general it will not have the promised injectivity properties. 
The problem originates in the number of connected components $V_{\eps^k}$, $\eps^k\in\Sigma_d^k$ 
is groving much faster than the number of connected components of $f^{-k}(U)$. 
To remove this problem we shall use decktransformations for $\phi$ 
to push together the sets $V_{\eps^k}$ and $V_{\wheps^k}$ whenever 
$\phi(V_{\eps^k})=\phi(V_{\wheps^k})$.

To fix the ideas let $\Gamma$ denote the group of decktransformations for the universal covering $\phi$, 
i.e~$\gamma\in\Gamma$, if and only if $\gamma$ is an automorphism of $\D$ with 
$\phi\circ\gamma =\phi$. 
\REFLEM{shiftedliftsofinverse}
Given $\epsk\in\Sigdk$ and  $\ga_1, \ldots , \ga_k\in \Gamma$ let 
$$
V = \ga_1\circ g_{\eps_1}\circ \ldots \circ \ga_k\circ g_{\eps_k}(\D)\quad\textrm{and}\quad
V' = \ga_2\circ g_{\eps_2}\circ \ldots \circ \ga_k\circ g_{\eps_k}(\D).
$$
Then the restrictions {\mapfromto {\phi} V {\phi(V)}} and {\mapfromto {\phi} {V'} {\phi(V')}} 
are universal coverings, $\phi(V)=W\Sm J_f$, $\phi(V')=W'\Sm J_f$, where 
$W, W'$ are connected components of $f^{-k}(U)$ and $f^{-(k-1)}(U)$ respectively and 
the restriction {\mapfromto f W {W'}} is a branched covering. 
\ENDLEM
\PROOF{}
The first statements is an easy induction proof, based on \propref{liftedinverses}, the details are left to the reader. 
The last statement follows from
$$
f\circ\phi\circ\ga_1\circ g_{\eps_1}\circ \ga_2 \circ\ldots \circ\ga_k\circ g_{\eps_k} 
= \phi\circ \ga_2\circ g_{\eps_2}\circ \ldots \ga_k\circ g_{\eps_k}.
$$
\ENDPROOF
Note also that 
$$
V\subset \ga_1\circ g_{\eps_1}\circ \ldots \circ \ga_{k-1}\circ g_{\eps_{(k-1)}}(\D).
$$

\REFPROP{goodlifts}
There exists a sequence of families ${\{\ga_l^{\epsk}\}}_{l=1}^k \subset \Gamma$, $\epsk\in\Sigdk$, $k\in\N$ 
such that the family of sets $V_{\epsk} := \ga_1^{\epsk}\circ g_{\eps_1}\circ \ldots \circ \ga_k^{\epsk}\circ  g_{\eps_k}(\D)$ 
and the sequence of families of decktransformations ${\{\ga_l^{\epsk}\}}_{l=1}^k$ satisfies the following three 
properties: 
\ENUM
\item\label{nested}
For all $k\geq 2$ and for all $\epsk$ : $V_{\epsk}\subset V_{\eps^{k-1}}$, where $\eps^{k-1}=\eps_1 \ldots \eps_{k-1}$.
\item\label{shiftinvariant}
For all $k\geq 2$, for all $\epsk$ and for all $l=2, \ldots , k$ : $\ga_l^{\epsk} = \ga_{l-1}^{\sigma(\epsk)}$.
\item\label{merged}
For all $k\geq 1$ and for all $\epsk,\whepsk$ : If $\phi(V_{\epsk}) = \phi(V_{\whepsk})$, then $V_{\epsk}= V_{\whepsk}$.
\ENDENUM
\ENDPROP
Remark that \itemref{merged} implies that there is a $1:1$ correspondence between connected components of 
$f^{-k}(U)$ and connected components of $\cup_{\epsk\in\Sigdk} V_{\epsk}$. And that \itemref{shiftinvariant} 
imples that this correspondence agrees with the dynamics, i.e.~$(f\circ\phi)(V_{\epsk}) = \phi(V_{\sigma(\epsk)})$.
\PROOF{}
The proof is by induction on $k$. For this it is convenient to let $\emptyset$ denote the empty tuple of length $0$ and define $\sigma(\eps^1)=\emptyset$. 
Also we shall then extextend the above properties \itemref{nested} and \itemref{shiftinvariant} to $k=1$ and 
property \itemref{merged} to $k=0$. 
We then define $V_\emptyset=\D$. This takes care of $k=0$. 
For $k=1$ we have already chosen the branches $g_j$ of the lifted inverse of $f$ so 
that $g_j(\D)=g_{j'}(\D)$ whenever $\phi(g_j(\D))=\phi(g_{j'}(\D))$. 
Thus we can simply take each $\ga_1^{\eps^1}=\id$. This then complies with all three properties.
For the inductive step suppose families ${\{\ga_l^{\epsn}\}}_{l=1}^n \subset \Gamma$, $\epsn\in\Sig_d^n$, 
$0\leq n<k$ satisfying the three properties have been constructed. 
For any $\epsk\in\Sigdk$ define $\ga_l^\epsk = \ga_{l-1}^{\sigma(\epsk)}$ for $1<l \leq k$. 
Moreover define $\whga_1^\epsk=\ga_1^{\eps^{k-1}}$ and 
$\whV_\epsk= \whga_1^\epsk\circ g_{\eps_1}(V_{\sigma(\epsk)})$ as preliminary candidates for 
$\ga_1^\epsk$ and $V^\epsk$. With this choice \itemref{shiftinvariant} is immediately satisfied and 
hence so is \itemref{nested}, because 
\begin{align*}
V_{\sigma(\epsk)} &=
\ga_1^{\sigma(\epsk)}\circ g_{\eps_2}\circ \ga_2^{\sigma(\epsk)}\circ g_{\eps_3}\circ\ldots\circ 
\ga_{k-1}^{\sigma(\epsk)}\circ g_{\eps_k}(\D)\\
&= \ga_2^{\epsk}\circ g_{\eps_2}\circ \ga_3^{\epsk}\circ g_{\eps_3}\circ\ldots\circ 
\ga_k^{\epsk}\circ g_{\eps_k}(\D)
\end{align*}
and by the induction hypothesis $V_{\sigma(\eps^{k-1})} \supset V_{\sigma(\epsk)}$, 
so that
$$
\whV_\epsk = \ga_1^{\eps^{k-1}}\circ g_{\eps_1}(V_{\sigma(\epsk)}) \subset 
 \ga_1^{\eps^{k-1}}\circ g_{\eps_1}(V_{\sigma(\eps^{k-1})}) = V_{\eps^{k-1}}.
 $$
 To complete the inductive step suppose $\phi(\whV_{\epsk}) = \phi(\whV_\whepsk)$. 
 Then by the above $\whV_\epsk \subset V_{\eps^{k-1}}$, 
 $\whV_\whepsk \subset V_{\wheps^{k-1}}$ so that
 $$
 \phi(V_{\eps^{k-1}}) = \phi(V_{\wheps^{k-1}}).
 $$
 And thus $V_{\eps^{k-1}} = V_{\wheps^{k-1}}$ by property \itemref{merged} applied to $\eps^{k-1}$ and 
 $\wheps^{k-1}$.
 That is $\phi(\whV_{\epsk}) = \phi(\whV_\whepsk)$ implies
 $$
 \whV_{\epsk}, \whV_{\whepsk} \subset V_{\eps^{k-1}} = V_{\wheps^{k-1}}.
 $$
 Define an equivalence relation $\sim$ on $\Sigdk$ by 
 $$
 \epsk\sim\whepsk \Leftrightarrow \phi(\whV_\epsk) = \phi(\whV_{\whepsk}).
 $$ 
 For each equivalence class of $\sim$ choose a prefered representative $\epsk$, 
 e.g. the one which is minimal with repsect to the lexicographic ordering, 
 and define $\ga_1^\epsk=\whga_1^\epsk$, $V_\epsk=\whV_\epsk$. 
 For any other element $\whepsk \in [\epsk]$ choose $\ga_1^\whepsk\in\Gamma$ so that 
 $$
 V_{\whepsk} = \ga_1^{\whepsk}\circ g_{\eps_1}(V_{\sigma(\whepsk)}) = V_{\epsk}.
 $$
 Then also property \itemref{merged} is satisfied.
 \ENDPROOF

We have now laid the grounds for the projection {\mapfromto \pi {\Sigd} {J_f}} of the Main Theorem:
Define the projection mapping {\mapfromto {\pi=\pi_f} {\Sigd} {J_f}} by
$$
\pi({(\eps_j)}_j) = \bigcap_{k=1}^\infty \overline{\phi(V_{\eps^k})}.
$$
Then by construction the map $\pi$ is continuous and semi-conjugates the shift $\sigma$ 
on $\Sigd$ to $f$ on $J_f$
$$
\pi\circ\sigma = f\circ\pi.
$$
The rest of the paper is devoted to proving that $\pi$ is as stated in the theorem: 
i.e. is surjective, is injective above any non-(pre)critical point $z$
and for any $z\in J_f$ satisfies
$$
\#\pi^{-1}(z) = \deg(f,z)\cdot\#\pi^{-1}(f(z)).
$$

Let us first address the issue of surjectivity. 
\REFPROP{surjectivity}
For any $k\in\N$ :
$$
f^{-k}(U\Sm J_f) = \bigcup_{\epsk\in\Sigdk} \phi(V_{\epsk})
$$
\ENDPROP
\PROOF{}
This is an elementary induction proof based on \propref{liftedinverses} and the observation that
for any $j$ and any decktransformation $\ga\in\Gamma$ : $\phi\circ\ga\circ g_j = \phi\circ g_j$. 
Combining the observation with \propref{liftedinverses} shows that the statement holds for $k=1$. 
Now suppose the statement holds for some $k$. Then by \propref{liftedinverses}:
$$
f^{-{(k+1)}}(U\Sm J_f) = \bigcup_{j=0}^{d-1}\bigcup_{\epsk\in\Sigdk} \phi\circ g_j(V_{\epsk}) 
= \bigcup_{\eps^{(k+1)}\in\Sig_d^{(k+1)}} \phi(V_{\eps^{(k+1)}}) 
$$
And thus the inductive step follows from the observation..
\ENDPROOF

\PROP The map {\mapfromto {\pi=\pi_f} {\Sigma_d} {J_f}} is surjective.
\ENDPROP
\PROOF{}
Let $z\in J_f$ be arbitrary and let $W_k$ denote the connected component of $f^{-k}(U)$ containing $z$. 
Then by \propref{surjectivity}, there exists a sequence ${(\eps^k)}_k$, $\eps^k\in\Sigma_d^k$ 
for each $k$ such that $\phi(V_{\epsk}) = W_k\Sm J_f$. 
By compactness of $\Sigma_d$ there exists at least one accumulation point $\veps\in\Sigma_d$ of
${(\eps^k)}_k$. That is for every $N\in\N$ there exists a $k>N$ such that 
$$
\eps_1, \ldots \eps_N = \epsk_1, \ldots \epsk_N.
$$
But then $\pi(\veps)\in W_k \subset W_N$ and thus $\pi(\veps)\in W_N$ for every $N$. 
That is $\pi(\veps) = z$, because $J_f$ is a Cantor set and thus $\diam(W_N)\to 0$ as 
$N\to \infty$. 
\ENDPROOF

\REFPROP{keyinjectivity}
If $V_{\epsk} = V_{\whepsk}$ and $\eps_1\not=\wheps_1$. Then the component $W$ 
of $f^{-k}(U)$ with $\phi(V_{\epsk}) = W\Sm J_f$ contains at least one critical point. 

Moreover if $\pi(\veps) = \pi(\whveps) = z$ and and $\eps_1\not=\wheps_1$. 
Then $z$ is a critical point for $f$.
\ENDPROP
\PROOF{}
Let $W' = f(W)$ then $\phi(V_{\sigma(\epsk)}) = \phi(V_{\sigma(\whepsk)}) = W'$ and 
thus $V_{\sigma(\epsk)} = V_{\sigma(\epsk)}$ by property \itemref{merged} of \propref{goodlifts}. 
Let $x$ be any point of the later set, then $\phi(g_{\eps_1})$ and $\phi(g_{\wheps_1})$ 
are two distinct preimages in $W$ of the point $\phi(x)\in W'$. 
Hence the degree of the restriction {\mapfromto f W {W'}} is at least $2$ and 
thus $W$ contains at least one critical point by the Riemann-Hurwitz formula.
This proves the first statement of the Lemma. 
The second is an immediate consequence of \propref{goodlifts} and the first statement:
\begin{align*}
\pi(\veps) &= \pi(\whveps)\\
\qquad &\Updownarrow \qquad\\
\bigcap_{k=1}^\infty \overline{\phi(V_{\epsk})} &= \bigcap_{k=1}^\infty \overline{\phi(V_{\whepsk})}\\
\qquad &\Updownarrow \qquad\\
\forall k\in\N : \phi(V_{\epsk}) &= \phi(V_{\whepsk})\\
\qquad &\Updownarrow \qquad\\
\forall k\in\N : V_{\epsk} &= V_{\whepsk}
\end{align*}
Thus if $\eps_1\not=\wheps_1$ and $W_k\Sm J_f = \phi(V_{\epsk})$. Then each $W_k$ contains a critical point, 
$W_{k+1}\subset\subset W_k$ for all $k$ and 
$$
z = \pi(\veps) = \cap_{k=1}^\infty \overline{\phi(V_{\eps^k})} = \cap_{k=1}^\infty W_k.
$$
Hence $z$ is a critical point.
\ENDPROOF

\COR 
Let $z\in J_f$ be any point whose orbit ${(f^n(z))}_{n\geq 0}$ does not contain a critical point. 
Then 
$$
\#\pi^{-1}(z) = 1.
$$
\ENDCOR
\PROOF{}
Suppose $\pi(\veps ) = \pi(\whveps) = z$. We shall show that $\veps = \whveps$. 
As a start $\eps_1=\wheps_1$ by \propref{keyinjectivity}. 
The Corollary now follows by induction since by the conjugacy property of $\pi$ 
$$
\pi(\sigma^n(\veps)) = \pi(\sigma^n(\whveps)) = f^n(z)
$$
and by assumption this point is not critical, so that $\eps_n=\wheps_n$ for all $n$ 
by \propref{keyinjectivity}.
\ENDPROOF

To shorten notation let us write $d_z = \deg(f,z)$ for any $z\in U'$.
\REFPROP{multiplicity}
For any $z\in J_f$ 
$$
\#\pi^{-1}(z) = d_z \#\pi^{-1}(f(z)).
$$
More precisely there are $d_z$ distinct numbers $j_1, \ldots j_{d_z}\in \{0, \ldots , d-1\}$ 
depending only on $z$ such that $\pi(\veps) = z$ if and only if  $\pi(\sigma(\veps)) = f(z)$ and 
$$
\eps_1 \in \{j_1, \ldots , j_{d_z}\}.
$$
\ENDPROP
\PROOF{}
Given $z\in J_f$ let $W_k$ denote the connected component of $f^{-k}(U)$ containing $z$ and 
let $W_k'$ denote the connected component of $f^{-(k-1)}(U)$ containing $f(z)$. 
Then the degree of the restrictions {\mapfromto f {W_k} {W_{k-1}'}} equals $d_z$ for $k$ sufficiently large, 
because $z$ is the only point in the nested intersection of the $W_k$. 
Fix any such $k_0$, let $k\geq k_0$ and let $V_{k-1}'$ denote the connected component of $\phi^{-1}(W_{k-1}'\Sm J_f)$ 
such that $V_{k-1}' = V_{\eps^{k-1}}$ for any $\eps^{k-1}$ with $\phi(V_{\eps^{k-1}}) = W_{k-1}'\Sm J_f$. 
Let $j_1, \ldots j_{d_z}\in \{0, \ldots , d-1\}$ be the $d_z$ values of $j$ 
for which $\phi(g_j(V_{k-1}')) = W_k\setminus J_f$ as provided by \propref{liftedinverses}. 
Then the index set $\{j_1, \ldots j_{d_z}\}$ does not depend on the value of $k\geq k_0$ by nestedness of 
the sets $V_{k-1}'$. 
Hence $\pi(\veps) = z$ if and only if $\eps_1\in\{j_1, \ldots j_{d_z}\}$ and $\pi(\sigma(\veps)) = f(z)$.
\ENDPROOF

\PROOF {\emph{(of \thmref{maintheorem})}} 
By the above Propositions $\pi$ is a continuous and surjective semiconjugacy.
In particular 
$$
\forall z\in J_f : \#\pi^{-1}(z) \geq 1.
$$
Since no critical point of $f$ is periodic and there are finitely many critical points counted with multiplicity, the total branching $\chi(z)$ along the orbit of an arbitrary point is uniformly bounded. In particular for any $z\in J_f$ there exists $N\in\N$ such 
that the orbit of $f^N(z)$ does not contain any critical point. Thus by \propref{keyinjectivity} 
$$
\#(\pi^{-1}(f^N(z))) = 1.
$$
Finally we have $\chi(z) = d_z\cdot d_{f(z)} \cdot \ldots \cdot d_{f^{N-1}(z)}$, so that 
$$
\forall z\in J_f : \#\pi^{-1}(z) = \chi(z).
$$
by induction on \propref{multiplicity}.
\ENDPROOF

Address:

Carsten Lunde Petersen, IMFUFA,
Roskilde University,
Postbox 260,
DK-4000 Roskilde,
Denmark.
e-mail: lunde@ruc.dk

\end{document}